\documentclass[12pt,twoside]{article}
\usepackage{amsmath,amsfonts,amstext,amsthm,array,multicol}
\usepackage[english]{babel}

\newcommand{\CC}{{\bf C}}
\newcommand{\cO}{{\cal O}}\newcommand{\cE}{{\cal E}}
\newcommand{\cD}{{\cal D}}
\newcommand{\gr}{{\rm gr}}

\newtheorem{theorem}{Theorem}[subsection]

\newtheorem{proposition}[theorem]{Proposition}
\newtheorem{lemma}[theorem]{Lemma}
\newtheorem{definition}[theorem]{Definition}

\title{Free divisors and duality for
${\cal D}$-modules \thanks{Partially supported by DGESIC-97-0723 and
HF-1998-0105. Second author partially supported by MSRI (Berkeley)}}
\author{F.J. Castro-Jim\'enez and J.M. Ucha}
\date{March 14, 2001}

\begin{document}
\maketitle

\abstract{The relationship between ${\cal D}$-modules and free divisors has
been studied in a general setting by L. Narv\'aez and F.J. Calder\'on. Using
the ideas of these works we prove in this article a duality formula between two
${\cal D}$-modules associated to a class of free divisors on ${\bf C}^n$ and we
give some applications.}

{\bf Math. Classification:} 32C38, 13N10,
14F40, 13P10.

\section{Free divisors}\label{free-divisors} Here we summarize some results of K. Saito
\cite{Saito}.

\setcounter{subsection}{1} Let us denote $X=\CC^n$. Denote by ${\cal O} = {\cal
O}_X$ the sheaf of holomorphic functions on $X$. Let $D\subset X$ be a divisor
and $x\in D$. Denote by $Der({\cal O}_x)$ the ${\cal O}_x$-module of ${{\bf
C}}$-derivations of ${\cal O}_x$ (the elements in $Der({\cal O}_x)$ are called
{\it vector fields}).

A vector field $\delta \in Der({\cal O}_x)$ is said to be {\it logarithmic}
w.r.t. $D$ if $\delta(f)=af$ for some $a\in {\cal O}_x$, where $f$ is a local
(reduced) equation of the germ $(D,x)\subset ({\bf C}^n,x)$. The ${\cal
O}_x$-module of logarithmic vector fields (or logarithmic derivations) is
denoted by $Der(\log D)_x$. This yields a $\cO$-module sheaf denoted by
$Der(\log D)$.

The divisor $D$ is said to be {\it free at the point} $x\in D$ if the ${\cal
O}_{x}$-module $Der(\log D)_x$ is free (and, in this case, of rank $n$). The
divisor $D$ is called {\it free} if it is free at each point $x\in D$.

Smooth divisors are free. A normal crossing divisor $D\equiv (x_1\cdots x_t=0)
\subset {\bf C}^n$ is free because we have $Der(\log D) = \bigoplus_{i=1}^t
{\cal O}_{{\bf C}^n}x_i\partial_i \oplus \bigoplus_{j=t+1}^n {\cal O}_{{\bf
C}^n}\partial_j$. By \cite{Saito} any reduced germ of plane curve $D\subset
{\bf C}^2$ is a free divisor.

Saito's criterium to test the freedom of a divisor $D$ at a point $p$ is:

\begin{lemma}({\rm \cite[(1.9)]{Saito}}) Let $\delta_i = \sum_{j=1}^n a_{ij} \partial_j$, $i=1,\ldots,n$
be a system of holomorphic vector fields at $p\in \CC^n$, such that: \\ i)
$[\delta_i,\delta_j] \in \sum_{k=1}^n {\cO}_p\delta_k$, for $i,j =1,\ldots,n$.
\\ ii) $det(a_{ij})=h$ defines a reduced hypersurface $D$. \\ Then, for
$D\equiv(h=0)$, $\delta_1,\ldots,\delta_n$ belong to $Der(\log D)_p$ and hence
$\{\delta_1,\ldots,\delta_n\}$ is a free basis of $Der(\log D)_p$.
\end{lemma}

\section{The logarithmic comparison theorem} \setcounter{subsection}{1}
Let $X$ be a complex manifold and $D\subset X$ a divisor. We have a canonical
inclusion $$ i_D : \Omega^\bullet(\log D) \rightarrow \Omega^\bullet(\star D)$$
where $\Omega^\bullet(\star D)$ is the meromorphic de Rham complex and
$\Omega^\bullet(\log D)$ is the de Rham  logarithmic complex, both w.r.t $D$. A
meromorphic form $\omega\in \Omega^p(\star D)$ is said to be {\it logarithmic}
if $fw\in \Omega^p$ and $df \wedge \omega \in \Omega^{p+1}$ for each local
equation $f$ of $D$.

A classical natural problem is to find the class of divisors $D\subset X$ for
which $i_D: \Omega^\bullet(\log D) \rightarrow \Omega^\bullet(\star D)$ is a
quasi-isomorphism (i.e. $i_D$ induces an isomorphism on cohomology).

By Grothendieck's comparison theorem we know that the complexes
$\Omega^\bullet(\star D)$ and ${\bf R}j_*({\bf C})$ are naturally
quasi-isomorphic, where $j:U=X\setminus D \rightarrow X$ is the natural
inclusion. So, if $i_D$ is a quasi-isomorphism we say that the logarithmic
comparison theorem holds for $D$ (or simply LCT holds for $D$).

\begin{definition}{\rm (\cite{CMN})} A divisor $D\subset X$ is locally quasi-homogeneous
if for all $q\in D$ there exist local coordinates $(V;x_1,\ldots,x_n)$ centered
at $q$ such that $D\cap V$ has a weighted homogeneous defining equation w.r.t.
$(x_1,\ldots,x_n)$. \end{definition}

Smooth divisors and normal crossing divisors are locally quasi-homogeneous. A
weighted homogeneous polynomial $f\in {\bf C}[x,y]$ defines a locally
quasi-homogeneous  divisor $D\equiv (f=0)\subset {\bf C}^2$.

Suppose $D\subset X$ is a locally quasi-homogeneous free divisor. The main
result of \cite{CMN} is that LCT holds for $D$, i.e. $$ i_D :
\Omega^\bullet(\log D) \rightarrow {\bf R}j_*({\bf C})$$ is a
quasi-isomorphism.

\section{Logarithmic $\cal D$-modules} Let us  denote  by ${\cal D} = {\cal D}_X$ the sheaf (of
rings) of linear differential operators with holomorphic coefficients on $X$.

A local section $P$ of ${\cal D}$ is a finite sum $$ P = \sum_\alpha a_\alpha
\partial^\alpha$$ where $\alpha=(\alpha_1,\ldots,\alpha_n)\in
{\bf N}^n$, $a_\alpha$ is a local section of ${\cal O}$ on some chart
$(U;x_1,\ldots,x_n)$ and
$\partial=(\partial_1,\ldots,\partial_n)=(\frac{\partial}{\partial x_1
},\ldots, \frac{\partial}{\partial x_n})$.

The sheaf $\cal D$ is naturally filtered by the order of the differential
operators. The associated graded ring ${\rm gr}({\cal D})$ is commutative. In
fact, we can identify ${\rm gr}({\cal D})$ with the sheaf ${\cal
O}[\xi_1,\ldots,\xi_n]$ of polynomials in the variables
$\xi=(\xi_1,\ldots,\xi_n)$ and with coefficients in ${\cal O}$.

Assume the operator $P=\sum a_\alpha \partial^\alpha$ has order $d$ (i.e.
$d=\max\{|\alpha|= \alpha_1+\cdots+\alpha_n\, \vert \, a_\alpha\not= 0 \}$)
then the {\it principal symbol} of $P$ is $$\sigma(P)= \sum_{|\alpha|=d}
a_\alpha \xi^\alpha \in {\cal O}[\xi].$$

For each left ideal $I$ in ${\cal D}$ the graded ideal associated to $I$ is the
ideal of ${\rm gr}({\cal D})$ generated by the set of principal symbol
$\sigma(P)$ for $P\in I$. This ideal is denoted by ${\rm gr}(I)$.

The {\it characteristic variety} of the ${\cal D}$-module $M=\frac{{\cal
D}}{I}$ is, by definition, the analytic sub-variety of the cotangent bundle
$T^*X$ defined by ${\cal O}_{T^*X}{\rm gr}(I)$. This characteristic variety if
denoted by $Ch(M)$. The cycle defined in $T^*X$ by the ideal ${\cal
O}_{T^*X}{\rm gr}(I)$ is denoted by $CCh(M)$ and it is called the {\it
characteristic cycle} of the ${\cal D}$-module $M$.

For any divisor $D\subset \CC^n$ the sheaf ${\cal O}[\star D]$ of meromorphic
functions with poles along $D$ is naturally a left coherent ${\cal D}$-module
(that follows from the results of Bernstein-Bj\"ork on the existence of the
$b$-function for each local equation $f$ of $D$, \cite{Bernstein},
\cite{Bjork}). Even more, Kashiwara proved that the dimension of $Ch({\cal
O}[\star D])$ is equal to $n$ (i.e. ${\cal O}[\star D]$ is holonomic,
\cite{Kas-II}).

In \cite{Cald} and \cite{Cald3} the author considers the (left) ideal $I^{\log
D }\subset \cD$  generated by the logarithmic vector fields $Der (\log D)$ (see
\ref{free-divisors}). We will denote simply $I^{\log}=I^{\log D}$ and
$M^{\log}$  the quotient $\cD/I^{\log}$ if no confusion is possible.

\subsection{Koszul free divisors}\label{Koszul-free}
Let us give the main result of F.J. Calder\'{o}n, \cite{Cald} (see also
\cite{Cald3}). Let $D\subset X$ be a divisor and $x\in D$.

\begin{definition}{\rm (\cite[Def. 4.1.1]{Cald3})}
The divisor $D$ is said to be Koszul free at the point $x\in D$ if it is free
at $x$ and there exists a basis $\{\delta_1,\ldots,\delta_n\}$ of $Der(\log
D)_x$ such that the sequence $\{\sigma(\delta_1),\ldots, \sigma(\delta_n)\}$ of
principal symbols is a regular sequence in the ring $\gr^F(\cD)$. The divisor
$D$ is Koszul free if it is Koszul free at any point of $D$.
\end{definition}

By \cite{Saito} and \cite[4.2.2.]{Cald3} any plane curve $D\subset \CC^2$ is a
Koszul free divisor. By \cite[Prop. 4.1.2]{Cald3} if $D$ is a Koszul free
divisor then $M^{\log}$ is holonomic and

\begin{theorem}\label{casi-iso-de-cald} {\rm (\cite[Th. 4.2.1]{Cald3})}
If $D$ is a Koszul free divisor then $\Omega^\bullet (\log D)$ and ${\bf
R}{\cal H}om_\cD(M^{\log},\cO)$ are naturally quasi-isomorphic.
\end{theorem}

\subsection{$\widetilde{M}^{\log}$}\label{Mtilde}
In \cite{Ucha-tesis} (see also \cite{Castro-Ucha-jsc}) L. Narv\'aez suggested
the study of the $\cD$-module $\widetilde{M}^{\log}$ defined as follows: Let us
denote by $\widetilde{I}^{\log}$ the left ideal of $\cD$ generated by the set
$\{ \delta + a\, \vert\, \delta\in I^{\log}{\mbox{ and }} \delta(f)=a f \}$.
Let us write $\widetilde{M}^{\log} = \cD/\widetilde{I}^{\log}$. There exists a
natural morphism $\phi_D : \widetilde{M}^{\log} \rightarrow \cO[\star D]$
defined by $\phi_D(\overline{P}) = P(1/f)$ where $\overline{P}$ denotes the
class of the operator $P\in \cD$ modulo $\widetilde{I}^{\log}$. The image of
$\phi_D$ is $\cD\frac{1}{f}$. As a natural question we ask for the class of $D$
such that the morphism $\phi_D$ is an isomorphism (see \ref{phi-D}).

\section{The duality theorem} \setcounter{subsection}{1}

Suppose here that the divisor $D\subset X$ is free, and let $f \in {\mathcal
O}$ be a local equation of $D$  and let $\{ \delta_1, ..., \delta_n \}$ be a
basis of the logarithmic derivations. We will use the following notation:
\begin{itemize}
\item $\delta_i(f) = m_i f$ for some $m_i\in \cO$.
\item $\delta_i = \sum_{k=1}^n{a_{ik} \partial_k}$ for some $a_{ik}\in \cO$.
\item $A = \left( \begin{array}{ccc}a_{11} & \cdots & a_{1n} \\
\vdots &  & \vdots \\ a_{n1} & \cdots & a_{nn} \end{array} \right)$
\item $[ \delta_i, \delta_j ] = \sum_{k=1}^n{\alpha^{ij}_k
\delta_k}$ for some $\alpha^{ij}_k \in \cO$.
\end{itemize}

\begin{lemma}{\label{det}} For any $i=1,...,n$ we have $$\delta_i (|A|) =
\sum_{k=1}^n{(\delta_i(a_{k1}), ..., \delta_i(a_{kn})) \left(
\begin{array}{c} A_{k1} \\ \vdots \\ A_{kn} \end{array}
\right)}$$ where $A_{kj}$ is the adjoint matrix of $a_{kj}$.
\end{lemma}
\begin{proof} From the very definition of the determinant
developed from the $k$-th row.
\end{proof}

The lemma above is true in fact for any derivation, not only for elements in
the basis.

\begin{lemma}{\label{alphas}} We have $$f (\alpha^{ij}_1, \ldots ,
\alpha^{ij}_n) = (\delta_i(a_{j1}) - \delta_j (a_{i1}), ..., \delta_i (a_{jn})
- \delta_j (a_{in})) Adj(A)^t.$$
\end{lemma}
\begin{proof}
It is only necessary to consider that $$[ \delta_i, \delta_j ] =
(\alpha^{ij}_1, ..., \alpha^{ij}_n) A \left(
\begin{array}{c}
\partial_1 \\ \vdots \\ \partial_n \end{array} \right) = $$
$$ = (\delta_i(a_{j1}) - \delta_j (a_{i1}), ..., \delta_i (a_{jn}) - \delta_j
(a_{in})) \left( \begin{array}{c}
\partial_1 \\ \vdots \\ \partial_n \end{array} \right).$$
\end{proof}

We consider the augmented Spencer logarithmic complex as in \cite[page
712]{Cald3}. We have $$ \cD\otimes_{\cO} \wedge^\bullet Der(\log D) \rightarrow
M^{\log} \rightarrow 0.$$

We say that a free divisor $D$ is {\it of Spencer type} if this complex is a
(locally) free resolution of $M^{\log}$ and this last $\cD$-module is
holonomic. By \cite[Prop. 4.1.3 ]{Cald3} if $D$ is Koszul free (in particular
if $D$ is a plane curve) then it is of Spencer type but the converse is not
true, see \cite[Remark 4.2.4]{Cald3} and section \ref{example-dim-3}.

The following proposition is a consequence of \cite[Th. 4.2.1]{Cald3}.
\begin{proposition}\label{casi-iso-de-cald-2}
If $D$ is of Spencer type then $Sol(M^{\log})$ is naturally quasi-isomorphic to
$\Omega^\bullet(\log D)$.
\end{proposition}

\begin{theorem}\label{duality-theorem}
Suppose $D$ is of Spencer type. Then $(M^{log})^* \simeq \widetilde{M}^{log}.$
\end{theorem}
\begin{proof} Using the Spencer logarithmic free resolution
of the holonomic $\cD$-module $M^{\log}$, we first compute a presentation of
the right $\cD$-module ${\cal E}:=Ext_\cD^n(M^{\log},\cD)$ and then we prove
that left $\cD$-module associated to $\cE$ is $\widetilde{M}^{\log}$.

The matrix of the $n$-th morphism in the resolution of $M^{log}$ (see
\cite[page 712]{Cald3}) has components of the form $$(- 1)^{i-1} \delta_i +
(-1)^{i}\sum_{l \neq i}{\alpha^{il}_l},$$ so it is enough to prove that
$$(-\delta_i + \sum_{k \neq i}{\alpha^{ik}_k})^* = \delta_i +
\sum_{k=1}^n{\partial_k (a_{ik})} + \sum_{k \neq i}{\alpha^{ik}_k} =  \delta_i
+ m_i.$$

In order to prove the last equality, we will show that $$m_i f = \delta_i (f) =
\delta_1 (|A|) =$$ $$= \sum_{k=1}^n {f
\partial_k (a_{ik})} + \sum_{k\neq i}{f \alpha^{ik}_k}.$$

Using \ref{alphas}, we obtain $$\sum_{k \neq i}{f \alpha^{ik}_k} = \sum_{k \neq
i}{(\delta_i(a_{k1}) - \delta_k (a_{i1}), ..., \delta_i (a_{kn}) - \delta_k
(a_{in})) \left(
\begin{array}{c} A_{k1} \\ \vdots \\ A_{kn}
\end{array} \right)} = $$ $$= \sum_{k = 1}^n {(\delta_i(a_{k1}), ..., \delta_i
(a_{kn})) \left(
\begin{array}{c} A_{k1} \\ \vdots \\ A_{kn} \end{array} \right)}
- \sum_{k =1}^n {(\delta_k (a_{i1}), ..., \delta_k (a_{in})) \left(
\begin{array}{c} A_{k1} \\ \vdots \\ A_{kn} \end{array} \right)}$$

So we have collected in the first sum precisely (see \ref{det}) $\delta_i
(|A|)$.  It remains to check that $$\sum_{k=1}^n {f
\partial_k (a_{ik})} = \sum_{k =1}^n {(\delta_k (a_{i1}), ..., \delta_k (a_{in}))
\left( \begin{array}{c} A_{k1} \\ \vdots \\ A_{kn}
\end{array} \right)}.$$

As $f = (a_{k1}, ..., a_{kn}) (A_{k1}, ..., A_{kn})^t$, we have $$\sum_{k=1}^n
{f \partial_k (a_{ik})} = \sum_{k=1}^n {
\partial_k (a_{ik})(a_{k1}, ..., a_{kn}) \left(
\begin{array}{c} A_{k1} \\ \vdots \\ A_{kn} \end{array}
\right)}= $$ $$= \sum_{k=1}^n {(\sum_{j=1}^n {a_{kj}
\partial_j (a_{i1})}, ..., \sum_{j=1}^n {a_{kj} \partial_j
(a_{in})}) \left(
\begin{array}{c} A_{k1} \\ \vdots \\ A_{kn} \end{array}
\right) } = $$ $$ = \sum_{k =1}^n {(\delta_k (a_{i1}), ..., \delta_k (a_{in}))
\left( \begin{array}{c} A_{k1} \\ \vdots
\\ A_{kn} \end{array} \right)}.$$
\end{proof}

\section{Some applications}

\subsection{LCT in dimension 2. Regularity of $M^{\log}$ and
$\widetilde{M}^{\log}$}\label{dim-2} Let $D\subset {\bf C}^2$ be a plane curve.

\begin{theorem} {\rm (\cite[Theorem 3.9]{4T})} The morphism $i_D :
\Omega^\bullet(\log D) \rightarrow \Omega^\bullet(\star D)$ is a
quasi-isomorphism if and only if $D$ is locally quasi-homogeneous.
\end{theorem}

\begin{proof} We show here how to read the original (topological) proof of
\cite{4T} to give a differential proof of ``only if" part. Part ``if" is a
consequence of \cite{CMN} because any plane curve is a free divisor
\cite{Saito}.

The problem is local. Suppose the local equation $f$ of $D$ is defined in a
small open neighbourhood such that the only singular point of $f=0$ is the
origin. Denote ${\cal O}[1/f]={\cal O}[\star D]$.

Let us consider (see \ref{Mtilde}) the natural surjective morphism $$ \phi_D :
\widetilde{M}^{\log} \rightarrow {\cal D}\frac{1}{f}\simeq {\cal O}[\star D]$$
where the last isomorphism follows by a result of Varchenko (i.e. the local
$b$-function $b_f(s)$ of $f$ verifies $b_f(-k)\neq 0$ for any integer $k\geq
2$, \cite{Va}). The kernel $K$ of $\phi_D$ is supported by the origin (because
$f$ is smooth outside $(0,0)$) and $CCh(\widetilde{M}^{\log}) = CCh(K) +
CCh({\cal O}[\star D])$. In particular $\widetilde{M}^{\log}$ and  $M^{\log} =
(\widetilde{M}^{\log})^*$ are regular holonomic (cf. \cite{MK2}) because as we
said before $D$ satisfies the hypothesis of theorem \ref{duality-theorem}.

Let us denote $Sol(M^{\log})={\bf R}{\cal H}om_{\cal D}(M^{\log},{\cal O})$ the
solution complex of $M^{\log}$.

Assume LCT holds for $D$. Then we have $$ DR({\cal O}[\star D]) \simeq
\Omega^\bullet(\star D) \simeq \Omega^\bullet(\log D) \simeq Sol(M^{\log})
\simeq DR((M^{\log})^*) \simeq DR(\widetilde{M}^{\log}).$$ Then both ${\cal
D}$-modules ${\cal O}[\star D]$ and $\widetilde{M}^{\log}$ have the same de
Rham complex and then the same characteristic cycle. In this case $K=0$ and
$\widetilde{M}^{\log} \simeq {\cal O}[\star D]$. Finally, by \cite[page
88]{Torrelli-tesis} (or by \cite[2.2.6]{Ucha-tesis}, see also
\cite{Castro-Ucha-jsc}) $f$ is weighted homogeneous in suitable coordinates.
That proves the ``only if" part of the theorem. \end{proof}

\subsection{On the comparison of $\widetilde{M}^{\log}$ and $\cO[\star
D]$}\label{phi-D}

In the previous section we proved (in dimension 2) that if
$\widetilde{M}^{\log} \simeq \cO[\star D]$ then  $f$ is weighted homogeneous
and the converse is also true (cf. \cite{Castro-Ucha-jsc}). So,
$\widetilde{M}^{\log} \simeq \cO[\star D]$ if and only if $f$ is weighted
homogeneous if and only if $LCT$ holds for $D$.

Now we return to dimension $n$.

\begin{theorem}
If $D\subset \CC^n$ is a free, locally quasi-homogeneous (l.q-h.) divisor then
$M^{\log}$ and $\widetilde{M}^{\log}$ are regular holonomic. Moreover
$\widetilde{M}^{\log D}$ and ${\cal O}[\star D]$ are naturally isomorphic.
\end{theorem}

\begin{proof} By \cite{Cald-Nar-1} $D$ is Koszul free and then $M^{\log}$ is
holonomic and the dual of $M^{\log}$ is $\widetilde{M}^{\log}$ (see
\ref{duality-theorem}). So, $\widetilde{M}^{\log}$ is also holonomic. It is
enough to prove that $\widetilde{M}^{\log}$ is regular.

To avoid confusion we will denote $\widetilde{M}^{\log D}$ to  emphasize the
divisor $D$.  In fact we will prove, by induction on $n$, that the natural
morphism $$\phi_D : \widetilde{M}^{\log D} \rightarrow \cO[\star D]$$ is an
isomorphism. We follows here the argument of \cite[4.3.]{Cald-Nar-2}. There is
nothing to prove in the case $n=1$. We note that in dimension 2 the result is
proved in \cite{Ucha-tesis} (see \ref{dim-2}). Suppose the result is true for
any free, l.q-h. divisor in dimension $\leq n-1$. Let $D\subset \CC^n$ be a
free, l.q-h. For any $x\in D$ there exists an open neighbourhood $U$ of $x$
such that for any $y\in U\cap D\setminus \{x\}$ the germ $(\CC^n,D,y)$ is
isomorphic to $(\CC^{n-1}\times \CC, D'\times \CC, (0,0))$, where $D'$ is a
free, l.q-h. divisor in $\CC^{n-1}$ (see \cite[prop. 2.4, lemma 2.2]{CMN}). So,
by induction hypothesis the morphism $\phi_{D'} : \widetilde{M}^{\log D'}
\rightarrow \cO[\star D']$ is an isomorphism. Then, by applying the functor
$\phi^*$ (where $\pi : \CC^n \rightarrow \CC^{n-1}$ is the projection), we have
that  for any $y\in U\cap D$, $y \neq x$, the morphism  $\phi_{D,y}$ is an
isomorphism between $\widetilde{M}^{\log D}_y$ and $\cO[\star D]_y$. We owe
this argument to L. Narv\'aez. So, the kernel of $\phi_D: \widetilde{M}^{\log
D} \rightarrow N^D$ is concentrated on a discrete set and it is regular
holonomic (here $N^D$ is the $\cal D$-module ${\cal D}\frac{1}{f}$, where $f$
is a local equation of $D$). As $N^D \subset {\cal O}[\star D]$ is regular
holonomic we deduce the regularity of $\widetilde{M}^{\log D}$. On the other
hand, by \cite{CMN} the logarithmic comparison theorem holds for $D$. So, by
using duality \ref{duality-theorem} and the natural quasi-isomorphism
$Sol(M^{\log D}) \rightarrow \Omega^\bullet(\log D)$ (\ref{casi-iso-de-cald}),
we deduce (as in \ref{dim-2}) that $DR(\widetilde{M}^{\log D})$ and $DR({\cal
O}[\star D])$ are naturally quasi-isomorphic and therefore, by Riemann-Hilbert
correspondence, $\widetilde{M}^{\log D}$ and ${\cal O}[\star D]$ are naturally
isomorphic, i.e. $\phi_D$ is an isomorphism. Thus we have concluded the
induction.
\end{proof}

\subsection{An example in dimension 3}\label{example-dim-3}

In \cite{4T} the authors give an example of a non Koszul free divisor --in
dimension 3-- for which LCT  holds. We will treat here, following the same
lines as in \cite{Castro-Ucha-jsc}, the case of the surface $D\subset \CC^3$
defined by $f= y(x^2+y)(x^2z+y)=0.$

The surface is free because computing the syzygies among $f, f_x, f_y, f_z$ we
obtain $$\begin{array}{c} (-3, \frac{1}{2}x, y, 0)
\\ (- x^2, 0, 0, x^2z+y) \\ (- xz - x, \frac{1}{2} x^2 +
\frac{1}{2} y, 0, xz^2 - xz), \end{array}$$ which produce the logarithmic
vector fields $$\begin{array}{ccl} \delta_1 & = & \frac{1}{2}x
\partial_x + y \partial_y \\ \delta_2 & = &  (x^2z+y)\partial_z \\ \delta_3 & = &
(\frac{1}{2} x^2 + \frac{1}{2} y) \partial_x + (xz^2 - xz)
\partial_z,
\end{array}$$
whose coefficients have a determinant equal to $1/2 f$.

This surface it is not Koszul-free because the set of the symbols (with respect
to the total order) $\sigma(\delta_1), \sigma(\delta_2), \sigma(\delta_3)$ do
not form a regular sequence. If we write $\sigma(\partial_x) = \xi,\
\sigma(\partial_y) = \eta,\ \sigma(\partial_3) = \zeta$. We have $y z \eta^2
\zeta + \frac{1}{4} \xi^2 \zeta \notin (\sigma(\delta_1), \sigma(\delta_2))$
but $ y z \eta^2 \zeta + \frac{1}{4} \xi^2 \zeta \sigma (\delta_3) \in
(\sigma(\delta_1), \sigma(\delta_2))$.

We compute a free resolution of $M^{log} = {\mathcal D}/{\mathcal D} (\delta_1,
\delta_2, \delta_3)$ using Gr\"obner basis. We obtain that the module of
syzygies $Syz(\delta_1, \delta_2, \delta_3)$ is generated by  ${\bf s}_{12},
{\bf s}_{13}, {\bf s}_{23}$ deduced from the commutators $[\delta_i, \delta_j
]$ where:
\begin{itemize}
\item $[\delta_1, \delta_2] = \delta_2$
\item $[\delta_1, \delta_3] = \frac{1}{2} \delta_3$
\item $[\delta_1, \delta_2] = (x z - x) \delta_2$.
\end{itemize}
The second module of syzygies (among the ${\bf s}_{ij}$) it is generated by
only one element ${\bf t}$, so we have finished the resolution. This element is
$${\bf t} = (t_1, t_2, t_3) = (-x^2z \partial_z + x z \partial_z - \frac{1}{2}
x^2 \partial_x - \frac{1}{2} y \partial_x - x z + x,\ x^2 z \partial_z + y
\partial_z, \ - y \partial_y - \frac{1}{2} x \partial_x + \frac{3}{2}),$$
precisely the one required in the Spencer logarithmic complex for $M^{\log}$ in
dimension 3, which is in fact a free resolution of $M^{\log}$. We can check
that in this case $M^{\log}$ is holonomic (use for example \cite{Macaulay2} or
\cite{Kan} ), so $D$ is of Spencer type and  we apply: \\ a) The theorem
\ref{duality-theorem} to obtain that $(M^{\log})^* \simeq
\widetilde{M}^{\log}$. \\ b) Proposition \ref{casi-iso-de-cald-2} to obtain
$Sol(M^{\log}) \simeq \Omega^\bullet(\log D)$.

Besides, the global $b$-function of $f$ is
$$(6s+5)(3s+2)(2s+1)(3s+4)(6s+7)(s+1)^3,$$ so we can assure that ${\mathcal
O}[\frac{1}{f}] \simeq \widetilde{M}^{\log}$, because $\widetilde{I}^{\log} =
Ann_{\mathcal D}(1/f)$.  We have used \cite{Macaulay2} and \cite{Kan} again to
compute the $b$-function and the annihilating ideal of $1/f$, that is to say,
the algorithms of \cite{Oaku}.

Finally, we have the following chain of quasi-isomorphisms $$ DR({\cal O}[\star
D]) \simeq \Omega^\bullet(\star D) \simeq \Omega^\bullet(\log D) \simeq
Sol(M^{\log}) \simeq DR((M^{\log})^*) \simeq DR(\widetilde{M}^{\log}).$$ and
the LCT holds for $D$.

\end{document}